\documentclass[12pt]{article}
\usepackage{amsmath,amsfonts}
\begin{document}

\baselineskip 20pt
\title{Radial limits and boundary uniqueness}

\author{Arthur A.~Danielyan}

\maketitle

\begin{abstract}

\noindent The paper
sheds a new light on the 
fundamental theorems
of complex analysis due to P. Fatou,
F. and M. Riesz, N. N. Lusin, I. I. Privalov, and A. Beurling.
 Only classical
tools
available at the times of Fatou are  used. The proofs are very simple and in some cases - almost trivial.

\end{abstract}

\begin{section}{Introduction and the main result.}

Let $U$ and $T$ be the open unit disk and the
unit circle in $\mathbb C$, recpetively.
For a function $f$ defined on $U$ we denote by
$f(e^{i \theta})$ the radial limit of $f$ at $e^{i \theta}$ if
the limit exists. The following classical theorems of P. Fatou \cite{fato}
(cf. \cite{colo}, Theorem 2.1) and of F. and M. Riesz \cite{rie}
(cf. \cite{colo}, Theorem 2.5) are among the most fundamental results of complex analysis.

\vspace{0.25 cm}

{\bf Theorem A (P. Fatou, 1906).}\ {\it Let $f$ be analytic and
bounded on $U$. Then for almost all $e^{i \theta}$ on $T$ the
radial limit $f(e^{i \theta})$ exists.}

\vspace{0.25 cm}

{\bf Theorem B (F. and M. Riesz, 1918).}\ {\it Let $f$ be analytic
and bounded on $U$ such that $f(e^{i \theta})=0$ on a set $E$ of
positive measure on $T$. Then $f$ is identically zero.}

\vspace{0.25 cm}

For univalent functions the analogous results are due to
A. Beurling \cite{beu} (cf. \cite{colo}, Theorem 3.5).

\vspace{0.25 cm}

{\bf Theorem C (A. Beurling, 1940).}\ {\it Let $f$ be univalent on
$U$. Then: (i) at every point $e^{i \theta}$ of $T$, except
possibly a set of zero (logarithmic) capacity, the
radial limit $f(e^{i \theta})$ exists; (ii) $f(e^{i
\theta})= \lim_{r \to 1} f(re^{i\theta})$ cannot be zero
on any positive capacity set on $T$.}

\vspace{0.25 cm}

Theorem A presents the property of almost everywhere existence of the radial limits of bounded analytic functions,
while Theorem B is
the boundary uniqueness property of the same functions. Both these properties (theorems)
have been presented and extended in
many books and countless papers.
Under the lights of the influential original works of Fatou and F. and M. Riesz, quite naturally, Theorem A
and Theorem B
have been universally regarded as two fundamental, but different properties of analytic functions\footnote{For instance the 
paper \cite{carl} by L. Carleson, emphasizing the difference between Theorem A and Theorem B, begins with the following sentences: 
 ``For a large number of classes $C$ of functions $f(z)$ regular in the unit circle, we have very complete knowledge 
 concerning the existence of a boundary function
  $$F({\theta})= \lim_{r \to 1} f(re^{i\theta}),$$ 
the classical result being that of Fatou. However, very little is known about the
properties of this boundary function $F(\theta)$, and in particular about the sets $E$ associated
with the class $C$, having the property that $f(z)$ vanishes identically if
$F(\theta) = 0$ on $E.$ ...
Our
whole knowledge in this direction seems to be contained in a classical result of F.
and M. Riesz: $E$ is a set of uniqueness for the class of bounded functions if and
only if it has positive Lebesgue measure."}.  
Undoubtedly the same persuasion has been the reason that the parts (i) and (ii) of Theorem C
have been regarded as different properties; obviously Beurling did not even suspect that one of the parts of Theorem C
might be a simple corollary of the other (for the proof of each part of Theorem C a special 
involved technic has been developed).

The  present paper sheds a new light on the
 theory and in fact changes
the mentioned viewpoint that Theorem A and Theorem B (or, part (i) and part (ii) of Theorem C) are presenting different properties. 
Instead, as we show, {\it the boundary uniqueness property is a direct corollary, or a simple particular case,
of the property of the a.e. existence of the radial limits.} 
In fact, except for trivial cases, {\it  the radial limits cannot be constant on a set $E$
just because the radial limits exist on a large subset of $E$ (or of $T$)},
 as the following main results (Theorem 1 and Corollary 1) of 
the paper show.

\vspace{0.25 cm}

{\bf Theorem 1.}\ {\it Let $f $ be univalent (respectively, zero free, bounded analytic)
on $U$ such that $f(e^{i \theta})=0$
on a subset $E$ of  $T$. Then $f$ generates a univalent (respectively,  bounded analytic) function $g$
on $U$ such that $g$ has no radial limit on $E$.}

\vspace{0.25 cm}

{\bf Remark 1.}\ If $f$ in Theorem 1 is zero free, bounded analytic, the relation between 
$f$ and $g$ is especially simple and given  explicitly by the equation
 $g(z)= e^{-i \log\log f(z)},$ as our below proof implies.

\vspace{0.25 cm}

Even in a trivial case of measure zero $E$, as we will see, 
Theorem 1 implies a theorem of Lusin \cite{lus}.
In contrast, the following corollary (of Theorem 1), in which $f$ is not required to be zero free on $U$,
 loses its meaning for $E$ of
Lebesgue measure zero.

\vspace{0.25 cm}

{\bf Corollary 1.}\ {\it Let $f $ be non-constant, bounded analytic
on $U$ such that $f(e^{i \theta})=0$ 
on a subset $E$ of $T$. Then $f$ generates a bounded analytic function $g$
on $U$ such that $g$ has no radial limit on $E$ except perhaps a subset of $E$ of
Lebesgue measure zero.}

\vspace{0.25 cm}

The hypothesis of Theorem 1 automatically implies that $f$ is non-constant.
But in Corollary 1 one needs to exclude the (trivial) case of the constant function.

Obviously Theorem 1 implies that the part (ii) of Theorem C is
a corollary of part (i) of Theorem C (nevertheless, the part (ii) itself is an important and frequently cited result).
Similarly, Corollary 1 derives Theorem B from Theorem A.  

Note that Theorem 1 (and Corollary 1) by implying both Theorem B and part (ii) of Theorem C, is giving the first
unified proof of the boundary uniqueness property for univalent and merely analytic functions.
Our approach greatly
simplifies 
the case of univalent functions and the new proof of part (ii) of Theorem C below should be compared 
with its previous proof; cf. the proof of part (ii) of Theorem C in \cite{colo}, pp. 61 - 64, or in  \cite{beu}.

In this paper we merely use
classical theorems (of times of Fatou) and modulo to them the proof of Theorem 1
is simple and elementary.

Theorem 1 
immediately implies also {\it the first elementary and very short proof} of the following theorem
 (see \cite{lus} or \cite{lupri}), which is a classical converse of Theorem A.

\vspace{0.25 cm}

{\bf Theorem D (N.N. Lusin, 1919).}\ {\it Let $E$ be a zero measure subset on $T$.
Then there exists a bounded analytic function $f$ on $U$ such that the radial limit $f(e^{i \theta})$
does not exist at each $e^{i \theta} \in E$.}

\vspace{0.25 cm}

We use the following classical result due to  Privalov \cite{Pri} (cf.  \cite{Priv}, p. 295, or \cite{Zyg}, p. 276).

\vspace{0.25 cm}

{\bf Theorem E (I.I. Privalov, 1919).}\ {\it Let $E$ be a zero measure subset on $T$.
Then there exists a zero free, bounded analytic, function $f$ on $U$ such that 
$f(z)$ tends to $0$ as $z$ approaches, in an arbitrary manner (in particlar, radially), any point of $E$.}

\vspace{0.25 cm}

Since the function  $f$
existing by this theorem is zero free, bounded analytic, and $f(e^{i \theta})=0$ on $E$,
it can serve as a function of hypotheses of Theorem 1, and thus, 
Theorem 1 readily implies Theorem D.
Under the light of our simple (below) proof of Theorem 1, Theorem D becomes nothing else but {\it an obvious  corollary} of Theorem E.
   Since also Privalov's
    proof of Theorem E is a simple  construction,
we arrive to the first elementary self-contained proof of Theorem D, which is in a sharp contrast to  its (very complex) original proof.

Lusin and Privalov have
  been collaborating for many years  
   in topics involving their Theorem D and Theorem E, and Theorem D even appears again
   in their well known joint paper \cite{lupri} (essentially with same original proof of 1919).  
    However they did not notice that Theorem D immediately
follows from  the elementary Theorem E.
(As a side effect of the present research, now one can finally include Theorem D with its new proof in the textbooks, and present it,
  along with Theorem A, in complex analysis graduate courses at the universities\footnote{Theorem A is proved in the standard textbooks on 
complex analysis, but, as a rule, none of them (say, W. Rudin's comprehensive  ``Real and Complex Analysis") even mentions Theorem D.
 This leaves the 
 reader wondering whether the conclusion of Theorem A is precise or not.
 An affirmative answer is provided by Theorem D and its 
  presence in the textbooks seems highly desirable.}).
  
 In this paper we also use two classical results due to C. Carath\'{e}odory (1913) and F. Riesz (1923), respectively.

\begin{subsection}{The theorems of Carath\'{e}odory and F. Riesz.}

 The prime end theorem of Carath\'{e}odory \cite{cara} (see for example, \cite{pom}, p. 30)
  in fact is the first major result on the boundary behavior of univalent functions.
It is on the extension of the Riemann mapping function to the boundary of a domain and
 has an especially simple formulation 
when the unit disc is mapped onto a Jordan domain. 
The general case uses the Carath\'{e}odory concept of a prime end of a simply connected domain $G$.
We assume that the reader is familiar with it as well as with the concept of a null-chain representing
the prime end (see for example \cite{pom}, pp. 29-30). However we will use just a special case of
Carath\'{e}odory's main theorem, which, as we will see, also can be derived from
the Riemann mapping theorem avoiding the prime ends (and null-chains) altogether.

\vspace{0.25 cm}

{\bf Theorem F (C. Carath\'{e}odory, 1913).}\ {\it Let $\varphi$ map the unit disc $U$ conformally onto a 
simply connected domain $G$. Then there is a bijective mapping $\hat \varphi$ of the circle $T$ onto the set of
all prime ends of $G$ such that, if $\zeta \in T$ and if $\{C_n\}$ is a null-chain representing the prime end
$\hat \varphi (\zeta)$, then $\{f^{-1}(C_n)\}$ is a null-chain that separates $0$ from $\zeta$ for large $n$.}

\vspace{0.25 cm}

If a Jordan arc in $U$ ends at a point $\zeta \in T$, then the image arc in $G$ ``approaches" (in the sense of prime ends)
 to the prime end
$\hat \varphi (\zeta)$, and vice versa.

Let $D$ be a simply connected ``double comb" 
domain in the $w$-plane obtained from the square
$$\{w=u+iv: 0<u<1, \ 0<v<1 \}$$ by taking off the line segments
 $l_{2n} =\{u+iv: u = \frac{1}{2n}, \ 0 \leq v \leq \frac{3}{4} \}$ and
$l_{2n+1} =\{u+iv: u = \frac{1}{2n+1}, \ \frac{1}{4} \leq v \leq 1 \}$ for all values of $n \  (n = 1, 2, ...)$.
Denote by $AB$ the closed set $\{iv: 0\leq v \leq 1 \}$ (the left side of the
original square). It contains no accessible boundary points of $D$.
In other words, there is no Jordan arc in $D$ ending at a point of $AB$ (to approach to $AB$, a Jordan
arc has to ``oscillate").
Note that $AB$ determines one prime end (of $D$), which we denote
 by $P$; more precisely,  $AB$ is the impression of $P$.
It is the only prime end of $D$ with an impression containing more than one point.

We have the following corollary of Theorem F.

\vspace{0.25 cm}

{\bf Proposition 1.}\ {\it Let $\varphi$  map $U$ conformally  onto $D$. 
Then there exists a point  $\xi$ on $T$ such that $\varphi$ has no limit as $z$ approaches $\xi$
along any Jordan arc $\gamma$, $\gamma \setminus \{\xi\} \subset U$, ending at $\xi$.}

\vspace{0.25 cm}

The existence of such $\xi$ immediately follows from Theorem F;
simply take as $\xi$ the point, which
corresponds  to the above mentioned prime end $P$ of $D$.

Let $\Gamma$ be a halfopen Jordan arc in $D$, oscillating and approaching to $AB$
asymptotically.
For instance, as such $\Gamma$, one can take the polygonal in $D$ joining the sequence of the
points $M_1(\frac{1}{2}, \frac{7}{8})$, $M_2(\frac{1}{3}, \frac{1}{8})$, $M_3(\frac{1}{4}, \frac{7}{8})$,
$M_4(\frac{1}{5}, \frac{1}{8})$,....
We may assume that   $\Gamma$ is
given by an equation $w=w(t), \ 0 \leq t <1$, where $w(t)$ is continuous on $[0,1)$ and  $w(0) \equiv M_1(\frac{1}{2}, \frac{7}{8})$ is the initial
point of $\Gamma$.
Let $\Gamma_1$ be a Jordan arc in $D$ having the same initial point $w(0)$ as $\Gamma$ and ending at an accessible
boundary point $w_1$ of $D$ ($\Gamma_1 \setminus \{w_1\} \subset D $), and such that $w(0)$
is the only common point of $\Gamma$ and $\Gamma_1$.

The following proposition is obvious.

\vspace{0.25 cm}

{\bf Proposition 2.}\ {\it The set $\Gamma \cup \Gamma_1$ divides the domain $D$ into
two domains $D_1$ and $D_2$ such that the boundaries of both $D_1$ and $D_2$ contain either all segments
$l_{2n}$ or all segments  $l_{2n+1}$ except finitely many of such segments.}

\vspace{0.25 cm}

Morera's theorem and the elementary (inner) uniqueness theorem immediately imply:

\pagebreak

{\bf Proposition 3.}\ {\it If $f$ is continuous in a domain $\Omega$ and analytic in $\Omega \setminus L$ where $L$ is a
line segment, then $f$ is analytic in $\Omega$. If in addition $f(z) = c$ on $L$, then $f$ is identically $c$ on $\Omega$.}

\vspace{0.25 cm}

For Proposition 1 we now present a simple, direct proof,
which avoids Theorem F and prime ends altogether, and only uses
the Riemann mapping theorem and
Theorem A. 

Let  $w=\varphi(z)$ be a conformal map of $U$ onto $D$ (as in Proposition 1).
Denote by $z=\psi(w)$ the inverse of $w=\varphi(z)$.
The image $\psi(\Gamma)$ of $\Gamma$
 is a halfopen Jordan arc in $U$ given by the equation $ w=\psi(w(t)), \ 0 \leq t <1$.
Note that $\psi(\Gamma)$ ends at a point $\xi \in T$, because otherwise
  $\psi(\Gamma)$ has to have two accumulating points $a$ and $b$ on $T$, and
   the function $\varphi$   
  cannot have 
  radial limits  
 on one of
  the two complementary to $a$ and $b$ open arcs of $T$,
  which contradicts to Theorem A. 
Now we show that $\xi \in T$ has the property formulated in Proposition 1.

  Let 
  $\gamma$, $\gamma \setminus \{\xi\} \subset U$,  
  be an arbitrary Jordan arc ending at $\xi$.
   If $\gamma$ and $\psi(\Gamma)$ share points
(other than $\xi$) at each neighborhood of $\xi$, then there is nothing to prove (because as
$\Gamma$, the curve $\varphi(\gamma)$ too would be oscillating and approaching to $AB$ in $D$).
Thus, by deleting some initial portion of $\gamma$ if necessary, we may assume that $\gamma$ and
$\psi(\Gamma)$ have no common point other than $\xi$. 
Let us join the initial points of $\gamma$ and
$\psi(\Gamma)$ by an arc $\delta \subset U$ such that $\delta$ has no other common point with $\gamma$
or with $\psi(\Gamma)$. 
The curve $\psi(\Gamma) \cup \delta \cup \gamma $ divides $U$ into two domains.
One of them, let denote it by $U_1$, has only one point of $T$, namely $\xi$, on its boundary.

Assume by contrary that $\varphi(z)$ has a limit equal to $q$ as $z$ approaches $\xi$ along $\gamma$.
This means that $q$ is an accessible boundary point of $D$ (and the Jordan arc  $\varphi(\gamma)$ ends at $q$).
The Jordan arc $\varphi(\delta)$ joins in $D$ the
initial points of $\Gamma$ and $\varphi(\gamma)$, and $\varphi(\delta)$ has no other common point
with them. Let us take $\varphi(\gamma) \cup \varphi(\delta)$ as
$\Gamma_1$ and apply Proposition 2; we conclude that there exists either a segment $l_{2n+1}$ or a segment
$l_{2n}$ lying on the boundary of the image $\varphi(U_1)$ of $U_1$. For briefness, denote this segment by $l$.

Because $\xi$ is the only boundary point of $U_1$ belonging to $T$, the conformal mapping $\psi(w)$ of
$\varphi(U_1)$ onto $U_1$ will be continuously extended to the set $l$ once we put $\psi(w) = \xi$ on $l$.
Now Proposition 3 implies that $\psi(w)$ is identically equal to $\xi$, which is impossible since $\psi$ is a univalent function.
This
contradiction completes the proof of Proposition 1.

\vspace{0.25 cm}

Next, we formulate a classical result on the radial limits of  Blaschke products.
By Theorem A, of course, the radial
limits of  a Blaschke product exist a.e. on $T$. 
In 1923
F. Riesz \cite{rie2} (cf. \cite{colo}, Theorem 2.11)
proved
 the following result.

\vspace{0.25 cm}

{\bf Theorem G (F. Riesz, 1923).}\ {\it A Blaschke product $B$ possesses radial limits of modulus $1$ for almost all $e^{i \theta}$ on $T$}

\vspace{0.25 cm}

 An elementary proof of this theorem
can be found in K. Hoffman's book\footnote {I am indebt to Don Marshall for
calling my attention to this proof.} (see \cite{hoff}, page bottom of 65 - top of 66). Without going into details, we describe the main steps of this proof. 
First, 
some simple integral estimates imply that for a (convergent) Blaschke product $B$ the sequence of partial products
 $B_n$ converge to $B$ in $H^2$ on the unit circle $T$. In particular, $B_n$ converges to $B$ in $L^2$ norm, and therefore,
a subsequence of $B_n$ converges pointwise a. e. on $T$ to $B$. Since each $|B_n|$ is identically $1$ on $T$, $|B|$ is $1$ a. e. on $T,$ as Theorem G claims.

We close this subsection with some remarks on our proof of Theorem 1. First of all we show that Proposition 1 easily implies Theorem 1. 
But, in fact, Proposition 1 is needed for Theorem 1 only for the case of univalent functions; 
for the case of just analytic functions we present yet another proof which instead of Proposition 1
merely uses an elementary argument.
Also, Theorem G immediately
reduces Corollary 1 to Theorem 1, and this is the only occasion when we use Theorem G in our proofs.

Thus, for the case of analytic functions, Theorem 1 is independent of Proposition 1,
and, in fact, Theorem 1 becomes almost trivial in this case (similarly, since $f$ provided by Theorem E is zero free, Theorem D is a trivial corollary of 
Theorem E, as mentioned above). 
For the case of univalent functions Theorem 1 follows almost trivially from Proposition 1.
And, Corollary 1 is almost trivial if we assume that Theorem G is granted.

{\it Thus in all cases
 the property of the existence of the radial limits almost trivially implies the boundary uniqueness property
(assuming at most Theorem G and Proposition 1 granted).}
(Based on below proofs, we hope it does make sense to use here the wording ``almost trivially", but the final judgement on this is reserved for the reader.)
Let us stress again that Proposition 1, for which we presented also a direct simple proof,
 is an obvious corollary of  Carath\'{e}odory's Theorem F of 1913, while Theorem G of F. Riesz 
 is of 1923 and has an elementary proof presented in \cite{hoff} as we noted above.

\end{subsection}

\end{section}

\begin{section}{Proofs.}

In the formulation of  Theorem 1 in general case the univalent function $f$ is not required to be bounded.
However, as it is well known, one can reduce this to the case of bounded $f$ using some 
 elementary mappings (cf. \cite{colo}, p. 57). Therefore, in the following proof with no loss of generality we assume that $f$ is bounded by $M$
 (on $U$) also for the case of univalent $f$.

\begin{subsection}{Proof of Theorem 1.}

Let $f$ be bounded, univalent (or analytic and zero free) on $U,$ and let
 $f(e^{i \theta})=0$ on $E.$
 (Since $f(e^{i \theta})=0$ on $E,$  of course, $f$ is zero free on $U$ also in
 case if $f$ is univalent.)
 We may assume
  $f$ is bounded by $1$. Then $f(z)=e^{h(z)}$, where $h$ is univalent (or analytic), $\Re h(z) < 0$ on $U$ and
 $h(e^{i \theta})=\infty$ on $E$. Let the univalent function $\varphi$ and the point $\xi \in T$ be as in Proposition 1, and let
  $\psi (\zeta)$ be a fractional-linear mapping of the left half plane
 onto the unit disk under which $\infty$ corresponds to $\xi$.
 Then the function $g(z)=\varphi(\psi(h(z)))$ does not have radial limits on $E$. 
 Next, $g$ is bounded analytic, and if $f$ is univalent, then with $h$ also $g$ is univalent. The proof is over.
 
\end{subsection}

\begin{subsection}{Elementary proof of Theorem 1 for the case of zero free,  bounded analytic $f$.}
 
We may assume $f$ is bounded by $1$. Then $f(z)=e^{h(z)}$, where $h$ is analytic, $\Re h(z) < 0$ on $U$ and
 $h(e^{i \theta})=\infty$ on $E$. We have an analytic $\log h(z)=\log |h(z)|+i \arg h(z)$  on $U$
  with $\pi/2<\arg h(z)<3\pi/2$.
 Then $g(z)= e^{-i \log h(z)}
 =e^{\arg h(z)}(\cos \log|h(z)| -i \sin \log |h(z)|)$ is analytic and bounded by $e^{3\pi/2}$.
On each radius ending on $E$
  the oscillation of  $g$ exceeds $e^{\pi/2}$.  The proof is over.

\end{subsection}

\begin{subsection}{Elementary proof of Lusin's Theorem D.}

 For a given zero measure set $E \subset T$ let $f$ be the function existing by Theorem E.
Repeating the previous paragraph for this $f$ proves Theorem D.

\end{subsection}

\begin{subsection}{Proof of Corollary 1.}

Let $f$ be as in Corollary 1.
Then  $f(z)=B(z)f_1(z)$ where $B$ is a Blashke product
 and $f_1$ is analytic, bounded, and zero free on $U$. By Theorem G we have $|B(e^{i \theta})|=1$ a.e. on $T$,
 and thus  $f_1(e^{i \theta})=0$ on some $E_1 \subset E$ such that $E \setminus E_1$ is of Lebesgue measure zero.
 By Theorem 1 for  $f_1$ there exists a bounded analytic function $g$ which does not have radial limits on $E_1$. The proof is over.
 
 \end{subsection}

\vspace{0.25 cm}

{\bf Remark 2.}\ Since $f_1$ is zero free, 
 $g(z)= e^{-i \log\log f_1(z)}$ as in Remark 1,
  and thus we have  $g(z)= e^{-i \log\log \frac{f(z)}{B(z)}}$ as an
  explicit equation which connects $f$ and $g$ of Corollary 1.
  
\vspace{0.4 cm}

{\bf Acknowledgment.} The author thanks Tom Carroll, Don Marshall, and Larry Zalcman for very helpful discussions
of a  preliminary version of this paper.

\end{section}

\begin{minipage}[t]{6.5cm}
Arthur A. Danielyan\\
Department of Mathematics and Statistics\\
University of South Florida\\
Tampa, Florida 33620\\
USA\\
{\small e-mail: adaniely@usf.edu}
\end{minipage}

\end{document}